\documentclass[12pt]{amsart}
\usepackage{amsmath,amscd,amssymb,amsfonts}
\def\ver{July 25, 2007, v.3}
\def\bA{{\bf A}}
\def\bC{{\bf C}}
\def\bN{{\bf N}}
\def\bP{{\bf P}}
\def\bQ{{\bf Q}}
\def\bR{{\bf R}}

\def\cD{{\mathcal D}}
\def\cG{{\mathcal G}}
\def\cI{{\mathcal I}}
\def\cJ{{\mathcal J}}
\def\cO{{\mathcal O}}
\def\fm{{\mathfrak m}}
\def\Re{{\rm Re}}
\def\ord{{\rm ord}\,}
\def\JN{{\rm JN}}
\def\RJN{{\rm RJN}}
\def\red{{\rm red}}
\def\codim{{\rm codim}\,}
\def\lct{{\rm lct}}
\def\rlct{{\rm rlct}}
\def\Sing{{\rm Sing}\,}

\begin{document}
\title{On real log canonical thresholds}
\author{Morihiko Saito}
\address{RIMS Kyoto University, Kyoto 606-8502 Japan}
\email{msaito@kurims.kyoto-u.ac.jp}
\date{\ver}
\begin{abstract}
We introduce real log canonical threshold and real jumping numbers
for real algebraic functions.
A real jumping number is a root of the $b$-function up to a sign if
its difference with the minimal one is less than 1. The real log
canonical threshold, which is the minimal real jumping number,
coincides up to a sign with the maximal pole of the distribution
defined by the complex power of the absolute value of the function.
However, this number may be greater than 1 if the codimension
of the real zero locus of the function is greater than 1.
So it does not necessarily coincide with the maximal root of the
$b$-function up to a sign, nor with the log canonical threshold of
the complexification.
In fact, the real jumping numbers can be even disjoint from the
non-integral jumping numbers of the complexification.
\end{abstract}
\maketitle

\centerline{\bf Introduction}

\bigskip\noindent
Let $f_{\bC}$ be a nonconstant holomorphic function on a complex
manifold $X_{\bC}$, and $\omega$ be a $C^{\infty}$ form of the
highest degree with compact support on $X_{\bC}$.
Then the integral $\int_{X_{\bC}}|f_{\bC}|^{2s}\omega$ is extended to
a meromorphic function in $s$ on the entire complex plane
(using a resolution of singularities [5] together with a partition of
unity, see [1], [2].)
Moreover, the largest pole of $\int_{X_{\bC}}|f_{\bC}|^{2s}\omega$
coincides up to a sign with the log canonical threshold of $f_{\bC}$
if $\omega$ is nonnegative and does not vanish on a point $x$ of
$D_{\bC}:=f_{\bC}^{-1}(0)$ where the log canonical threshold of
$(f_{\bC},x)$ attains the minimal.
(This follows from the definition by using a resolution of
singularities, see [7].)

Let $f$ be a nonconstant real algebraic function on a real algebraic
manifold $X_{\bR}$, and $\omega$ be a $C^{\infty}$ form of the
highest degree with compact support on $X_{\bR}$ such that the open
subset $\{x\in X_{\bR}\,|\,\omega(x)\ne 0\}$ is oriented and
$\omega(x)$ is positive on this subset.
Then $\int_{X_{\bR}}|f|^s\omega$ is similarly extended
to a meromorphic function in $s$ on the entire complex plane.
But the largest pole of $\int_{X_{\bR}}|f|^s\omega$
does not necessarily coincide up to a sign with the log canonical
threshold of the complexification $f_{\bC}:X_{\bC}\to\bC$ of
$f:X_{\bR}\to\bR$, see Corollary 2 and Theorem~1 below.

Let $\cO_{X_{\bR}}$ denote the sheaf of real analytic functions on
$X_{\bR}$.
We define the {\it real multiplier ideals}
$\cJ(X_{\bR},f^{\alpha})\subset\cO_{X_{\bR}}$ for $\alpha\in\bQ_{>0}$
by the local integrability of $|g|/|f|^{\alpha}$ for
$g\in\cO_{X_{\bR}}$.
(Here coherence of $\cJ(X_{\bR},f^{\alpha})$ is unclear.)
We have $\cJ(X_{\bR},f^{\alpha})=\cO_{X_{\bR}}$ for $0<\alpha\ll 1$,
but not necessarily
$$
f\cJ(X_{\bR},f^{\alpha})=\cJ(X_{\bR},f^{\alpha+1})\,\,\,\hbox{for}
\,\,\,\alpha>0,
$$
unless $f$ is of ordinary type.
Here we say that $f$ is {\it of ordinary type} if $\codim D_{\bR}=1$
where $D_{\bR}=f^{-1}(0)\subset X_{\bR}$,
and {\it of exceptional type} otherwise.
Note that the above equality always holds in the complex case.

By Hironaka [5], there is a resolution of singularities as real
algebraic manifolds $\pi_{\bR}:X'_{\bR}\to X_{\bR}$ which is a
composition of blowing-ups along smooth centers over $\bR$
and such that $\pi^*f$ and $\pi^*dx_1\cdots dx_n$ are locally of
the form $u\prod_{i=1}^rx'_i{}^{m_i}$ and $u'\prod_{i=1}^rx'_i
{}^{a_i}dx'_1\cdots dx'_n$respectively,
where $m_i\ge 1$ for $i\in[1,r]$.
Here $x_1,\dots, x_n$ and $x'_1,\dots,x'_n$ are local coordinates of
$X_{\bR}$ and $X'_{\bR}$ respectively, and $u,u'$ are nowhere
vanishing.
So $\pi^*f$ defines a divisor with normal crossings $D'_{\bR}=
\sum_{j\in J_{\bR}}m_jD'_{j,\bR}$, and we may assume that each
$D'_{j,\bR}$ is smooth by loc.~cit.
Let $a_j$ be the multiplicity of the Jacobian of $\pi_{\bR}$ along
$D'_{j,\bR}$.
Note that $m_j$ and $a_j$ are given by the above $m_i$ and $a_i$
respectively if $D'_{j,\bR}$ is locally defined by $y'_i=0$.

\medskip\noindent
{\bf Proposition~1.} {\it
For $g\in\cO_{X_{\bR},x}$ we have}
$$\hbox{$g\in\cJ(X_{\bR},f^{\alpha})_x\iff \pi^*gdx_1\cdots dx_n\in
\big(\pi_*\Omega^n_{X'_{\bR}}(-\sum_j[\alpha m_j]D'_{j,\bR})\big)_x$}.
$$

\medskip
However, $\pi_*\Omega^n_{X'_{\bR}}(-\sum_j[\alpha m_j]D'_{j,\bR})$
may be larger than $\cJ(X_{\bR},f^{\alpha})\Omega^n_{X_{\bR}}$
in general (even for $0<\alpha\ll 1)$, and coherence of
these sheaves are unclear.
By Proposition~1 there are increasing rational numbers
$0<\alpha_1<\alpha_2<\cdots$ such that
$$
\cJ(X_{\bR},f^{\alpha_j})=\cJ(X_{\bR},f^{\alpha})\supsetneq
\cJ(X_{\bR},f^{\alpha_{j+1}})\,\,\,\hbox{if}\,\,\,
\alpha_j\le\alpha<\alpha_{j+1}\,\,\,(j\ge 1),
$$
and $\cO_{X_{\bR}}=\cJ(X_{\bR},f^{\alpha})\supsetneq\cJ(X_{\bR},
f^{\alpha_1})$ if $0<\alpha<\alpha_1$.
These numbers $\alpha_j$ are called {\it real jumping numbers} of $f$.
(Here we add ``real'' since the complexification $f_{\bC}$ of $f$ can
be identified with $f$ in case $f\in\bR[x]\subset\bC[x]$.)
The minimal real jumping number $\alpha_1$ is called the {\it real
log canonical threshold}, and is denoted by $\rlct(f)$.
This is the smallest number such that $|f|^{-\alpha}$ is not locally
integrable on $X_{\bR}$.
It may be strictly greater than 1 in case of exceptional type,
see Theorem~1 below.
We define the graded pieces by
$$
\cG(X_{\bR},f^{\alpha})=
\cJ(X_{\bR},f^{\alpha-\varepsilon})/\cJ(X_{\bR},f^{\alpha})\quad
\hbox{for}\,\,\,0<\varepsilon\ll 1,
$$
so that $\alpha$ is a real jumping number of $f$ if and only if
$\cG(X_{\bR},f^{\alpha})\ne 0$.
Proposition~1 implies

\medskip\noindent
{\bf Corollary~1.} {\it We have}
$$
\rlct(f)=\min_{j\in J_{\bR}}\Big\{\frac{a_j+1}{m_j}\Big\}.
$$

\medskip
A similar assertion holds for the log canonical threshold
$\lct(f_{\bC})$ by applying the same argument to the resolution of
singularities of the complexification $f_{\bC}$,
and $-\lct(f_{\bC})$ coincides with the largest root of
$b_{f_{\bC}}(s)$, see [7].
Let $-p(f,\omega)$ denote the maximal pole of
$\int_{X_{\bR}}|f|^s\omega$.
Then

\medskip\noindent
{\bf Corollary~2.} {\it
We have in general
$$
p(f,\omega)\ge\rlct(f)\ge\lct(f_{\bC}),
$$
and $p(f,\omega)=\rlct(f)$ if $\omega(x)\ne 0$ for some $x\in X_{\bR}$
such that $\cG(X_{\bR},f^{\alpha})_x\ne 0$ with $\alpha=\rlct(f)$.}

\medskip
For the corresponding assertion in the complex case, see [7].
The relation with the complexification is quite complicated
as is shown by the following.

\medskip\noindent
{\bf Theorem~1.} {\it
There are cases where $\rlct(f)>\lct(f_{\bC})$, and even
$\rlct(f)>1$ in case of exceptional type.
Moreover the real jumping numbers of $f$ can be disjoint from the
non-integral jumping numbers of $f_{\bC}$ even in the case $f_{\bC}$
has only an isolated singularity at a real point
$x\in X_{\bR}\subset X_{\bC}$.}

\medskip
This kind of phenomena may happen in case $f$ has an isolated zero of
simple type, see (3.3).
Let $b_f(s)$ be the $b$-function of $f$ which is by definition
the least common multiple of the local $b$-functions $b_{f,x}(s)$
for $x\in X_{\bR}$.
Note that $b_{f,x}(s)$ coincides with the local $b$-function
$b_{f_{\bC},x}(s)$ of $f_{\bC}$,
since $b_{f_{\bC},x}(s)\in\bQ[s]$ by Kashiwara~[6].
So $b_f(s)=b_{f_{\bC}}(s)$ in case $\Sing f\subset X_{\bR}$.

\medskip\noindent
{\bf Theorem~2.} {\it
Any real jumping number of $f$ which is smaller than $\rlct(f)+1$ is
a root of $b_f(-s)$.}

\medskip
For the corresponding assertion in the complex case, see [4].
It seems that the case of an ideal generated by $f_1,\dots,f_r$
is reduced to the case $r=1$ by considering
$f=\sum_{i=1}^rf_i^2$ in the real case.

This note is written to answer questions of Professor S.~Watanabe
which are closely related to problems in the theory of learning
machines (see e.g. [10]).
I would like to thank him for interesting questions.

In Section 1 we recall some facts from the theory of
resolutions of singularities due to Hironaka [5].
In Section 2 we prove Proposition~1 and Theorem~2.
In Section 3 we prove Theorem~1 by constructing examples.

\bigskip\bigskip
\centerline{\bf 1. Resolution of singularities}

\bigskip\noindent
In this section we recall some facts from the theory of
resolutions of singularities due to Hironaka [5].

\medskip\noindent
{\bf 1.1.~Analytic spaces associated to $\bR$-schemes.}
Let $X$ be a scheme of finite type over $\bR$.
We denote the associated real analytic space by $X_{\bR}$.
The underlying topological space of $X_{\bR}$ is the set of
$\bR$-valued points $X(\bR)$ with the classical topology.
The sheaf of real analytic functions on $X_{\bR}$ is defined by
taking local embeddings of $X$ into affine spaces and dividing
the sheaf of real analytic functions on the affine spaces
by the corresponding ideal.

We define $X_{\bC}$ similarly for a scheme $X$ of finite type over
$\bC$.
In case $X$ is a scheme of finite type over $\bR$, $X_{\bC}$ means
the complex algebraic variety associated to the base change of
$X$ by $\bR\to\bC$.
So the underlying topological space of $X_{\bC}$ coincides with
$X(\bC)$.

\medskip\noindent
{\bf 1.2.~Hironaka's resolution of singularities.}
Let $X$ be a smooth scheme over $\bR$, and $D$ an effective divisor
on $D$.
By Hironaka [5] we have a resolution of singularities
$\pi:(X',D')\to(X,D)$ which is a composition of blowing-ups along
smooth centers defined over $\bR$ and such that $D'$ is a divisor
with normal crossings which is locally defined by {\it algebraic}
local coordinates defined over $\bR$, see loc.~cit., Cor.~3 in p.~146
and also Def.~2 in p.~141.
(Note that the last condition implies that the irreducible components
$D'_j$ of $D'\,(j\in J)$ are smooth over $\bR$ by taking a point of
$\Sing D'_j$).

This induces a resolution of singularities
$\pi_{\bR}:(X'_{\bR},D'_{\bR})\to(X_{\bR},D_{\bR})$ as in
Introduction, and
$$
J_{\bR}=\{j\in J\,|\,D'_j(\bR)\ne\emptyset\}.
$$
Note that if a smooth center $C$ of a blow-up has a real point $x$,
then $C$ is defined locally by using local algebraic coordinates
over $\bR$, and hence $C_{\bR}$ is a smooth subvariety.

\bigskip\bigskip
\centerline{\bf 2. Proofs of Proposition~1 and Theorem~2}

\bigskip\noindent
In this section we prove Proposition~1 and Theorem~2.

\medskip\noindent
{\bf 2.1.~Proof of Proposition~1.}
With the notation of Introduction, we have locally
$$
\hbox{$\pi^*g f^{-\alpha}dx_1\cdots dx_n=
v\prod_{i=1}^rx'_i{}^{a_i+b_i-\alpha m_i}dx'_1\cdots dx'_n$},
$$
if $\pi^*g=u''\prod_i x'_i{}^{b_i}$ locally, where $v,u''$
are nondivisible by $x'_i\,(1\le i\le r)$.
For $\gamma,c>0$, we have
$$
\int_0^{c}x{}^{\gamma-1}dx=\frac{c^{\gamma}}{\gamma},
$$
where $x$ means $x'_i$.
Moreover, we have for $\beta=\alpha m_i$ and
$p=a_i+b_i$
$$
p\ge[\beta]\iff p>\beta-1.
\leqno(2.1.1)
$$
So the implication $\Leftarrow$ in Proposition~1 follows.
For the converse, assume the right-hand side does not hold.
Then the left-hand side does not hold by restricting to a
neighborhood of a sufficiently general point of $D'_{j,\bR}$
which is defined locally by $x'_i=0$ and such that
$a_i+b_i-\alpha m_i\le-1$ (using positivity).
So the assertion follows.

\medskip\noindent
{\bf 2.2.~Proof of Corollary~1.}
By definition the minimal real jumping number is the smallest number
$\alpha$ such that $1\notin\cJ(X_{\bR},f^{\alpha})$, i.e.
$|f|^{-\alpha}$ is not locally integrable on $X_{\bR}$.
By Proposition~1, this condition is equivalent to that
$a_j<[\alpha m_j]$ (i.e. $a_j\le \alpha m_j-1$, see (2.1.1)) for some
$j\in J_{\bR}$.
So the assertion follows.

\medskip\noindent
{\bf 2.3.~Proof of Corollary~2.}
We take a resolution of singularities as in (1.2).
This gives a resolution of singularities of the complexification.
We define similarly $a_j,m_j$ for any irreducible components $D'_j$
of $D'\,(j\in J)$, and we have as in [7]
$$
\lct(f_{\bC})=\min_{i\in J}\Big\{\frac{a_j+1}{m_j}\Big\}.
$$
So the last inequality follows.
Since $\rlct(f)$ is the smallest number $\alpha$ such that
$|f|^{-\alpha}$ is not locally integrable on $X_{\bR}$,
the first inequality and the last assertion follow.

\medskip\noindent
{\bf 2.4.~Proof of Theorem~2.}
Let $f_+(x)=f(x)$ if $f(x)>0$ and $f_+(x)=0$ otherwise.
Set $f_-=(-f)_+$.
Since $|f|^s=(f_+)^s+(f_-)^s$, we consider
$$
I(\omega,s)=\int_{X_{\bR}}(f_+)^s\omega,
$$
where $\omega$ is a $C^{\infty}$ form of the highest degree
whose support is compact and is contained in a sufficiently
small open subset $U_{\bR}$ of $X_{\bR}$ with local coordinates
$x_1,\dots,x_n$ giving an orientation of $U_{\bR}$.
Then $I(\omega,s)$ is a holomorphic function on
$\{s\in\bC\,|\,\Re\,s>0\}$, and it is extended to a meromorphic
function on the entire complex plane using a resolution of
singularities, see [1], [2].

Let $x$ be a point of $D_{\bR}:=f^{-1}(0)\subset X_{\bR}$,
and $b_f(s)$ be the $b$-function of $f$ at $x$.
We assume that $U_{\bR}$ is a sufficiently small open neighborhood
of $x$ in $X_{\bR}$ so that we have the relation
$$
b_{f}(s)f^s=Pf^{s+1}\,\,\,\text{in}\,\,\,(\cO_{U_{\bR}}
\hbox{$[\frac{1}{f}])$}
[s],\,\,\,\text{where}\,\,\, P\in \Gamma(U_{\bR},\cD_{U_{\bR}}[s]).
\leqno(2.4.1)
$$
Here $P$ is replaced by $-P$ if $f_+$ is replaced by $f_-$
(and $f$ by $-f$).
Note that (2.4.1) holds in $\cO_{U_{\bR}}[\frac{1}{f}]$
when $s$ is specialized to any complex number.

Let $*$ be the involution of $\cD_{U_{\bR}}$ such that
$g^*=g$ for $g\in\cO_{U_{\bR}}$, $(\partial/\partial x_i)^*=
-\partial/\partial x_i$, and $(Q_1Q_2)^*=Q_2^*Q_1^*$ for
$Q_1,Q_2\in\cD_{U_{\bR}}$, fixing the local coordinates
$x_1,\dots,x_n$ on $U_{\bR}$.
This gives a right $\cD_{U_{\bR}}$-module structure
on $\Omega_{U_{\bR}}^n$ using a basis $dx_1\wedge\cdots\wedge dx_n$.
Write $P=\sum_jP_js^j$ with $P_j\in\cD_{\bR}$, and set
$P^*=\sum_jP^*_js^j$.
Let $r=\max\{\ord P_j\}$.
Then, for any complex number $s$ with $\Re\,s>r$, we have by (2.4.1)
together with integration by parts
$$
b_{f}(s)I(\omega,s)=\int_{U_{\bR}}b_{f}(s)(f_+)^s\omega=
\int_{U_{\bR}}(f_+)^{s+1}(P^*\omega)=\sum_jI(P_j^*\omega,s+1)s^j,
\leqno(2.4.2)
$$
since $\prod_i(\partial/\partial x_i)^{\nu_i}(f_+)^s$ is a
continuous function on $U_{\bR}$ if $\Re\,s>\sum_i\nu_i$.
Here $P_j^*\omega$ is defined by trivializing $\Omega_{U_{\bR}}^n$
by $dx_1\wedge\cdots\wedge dx_n$, and it may be written as
$\omega P_j$ using the right $\cD$-module structure explained above.
By analytic continuation, (2.4.2) holds as meromorphic functions in
$s$ on the entire complex plane.

Let $\alpha$ be a real jumping number of $f$ which is smaller than
$\rlct(f)+1$.
Assume that the above $x$ belongs to the support of
$\cG(X_{\bR},f^{\alpha})$, and $\omega(x)\ne 0$.
There is $g\in\Gamma(U_{\bR},\cO_{U_{\bR}})$ such that
$g\in\cJ(U_{\bR},f^{\alpha-\varepsilon})_x$ for $\varepsilon>0$ and
$g\notin\cJ(U_{\bR},f^{\alpha})_x$
(shrinking $U_{\bR}$ if necessary).
Then $I(g\omega,s)$ is a holomorphic function in $s$ for
$\Re\,s>-\alpha$ using a resolution of singularities as in (2.1),
and
$$
I(g\omega,s)\to +\infty\,\,\,\hbox{as}\,\,\,
s\to -\alpha,
$$
(replacing $f_+$ with $f_-$ if necessary).
On the other hand, the $I(P_j^*(g\omega),s+1)$ are holomorphic
functions in $s$ for $\Re\,s+1>-\rlct(f)$.
Thus, replacing $\omega$ with $g\omega$ in (2.4.2), we get
$b_f(-\alpha)=0$ since $-\alpha+1>-\rlct(f)$.
So the assertion follows.

\medskip\noindent
{\it Remark.}
This argument shows that the order of pole of $I(\omega,s)$ at
$-\rlct(f)$ is at most the multiplicity of $-\rlct(f)$ as a root
of $b_f(s)$.

\medskip\noindent
{\bf 2.5.~$b$-Function of the complexification.}
For $f\in\bR\{\!\{x\}\!\}$, the $b$-function $b_f(s)$ of
$f$ coincides with the $b$-function $b_{f_{\bC}}(s)$ of the
complexification $f_{\bC}$ (which is identified with $f$ by
$\bR\{\!\{x\}\!\}\subset\bC\{\!\{x\}\!\}$), since
$b_{f_{\bC}}(s)\in\bQ[s]$ by Kashiwara~[6].

Indeed, if there is
$P=\sum_{\nu,\mu,k}a_{\nu,\mu,k}x^{\nu}\partial^{\mu}s^k$
with $a_{\nu,\mu,k}\in\bC$ and satisfying
$$
b_{f_{\bC}}(s)f^s=Pf^{s+1},
$$
then the same equation holds with $P$ replaced by
$\sum_{\nu,\mu,k}(\Re\,a_{\nu,\mu,k})x^{\nu}\partial^{\mu}s^k$.

\medskip\noindent
{\bf 2.6.~Case of ideals.} For an ideal $\cI$ generated
by $f_1,\dots,f_r$, we may define the multiplier ideals
$\cJ(X_{\bR},\cI^{\alpha})$ by local integrability of
$$
\hbox{$|g|/(\sum_i|f_i|)^{\alpha}$}.
$$
However, this is calculated by
$\cJ(X_{\bR},f^{\alpha/2})$ with $f=\sum_if_i^2$, using
$$
\hbox{$\sum_i|f_i|^2\le(\sum_i|f_i|)^2\le r\sum_i|f_i|^2$}.
$$

\bigskip\bigskip
\centerline{\bf 3. Proof of Theorem~1}

\bigskip\noindent
In this section we prove Theorem~1 by constructing examples.

\medskip\noindent
{\bf 3.1.~Definition.}
We say that $f$ is {\it of ordinary type} if $\codim D_{\bR}=1$,
and {\it of exceptional type} otherwise.
Here $D_{\bR}=f^{-1}(0)\subset X_{\bR}$.

Write $f=\sum_{k\ge d}f_k\in\bR\{\!\{x_1\dots,x_n\}\!\}$ with $f_k$
homogeneous of degree $k$ and $f_d\ne 0$.
We say that $f$ has an {\it isolated zero of simple type} if the
equation $f_d=0$ has no solution in $\bR^n\setminus\{0\}$
(e.g. if $f_d=\sum_{i=1}^nx_i^d$ with $d$ even).

\medskip\noindent
{\bf 3.2.~Remarks.} (i)
The function $f$ is of ordinary type if and only if the reduced
complex zero locus $(D_{\bC})_{\red}$ has a smooth real point.
Note that
$$
\dim_{\bR}(D_{\bR}\cap\Sing (D_{\bC})_{\red})<n-1,
$$
since $\Sing (D_{\bC})_{\red}$ is defined over $\bR$ and has
dimension $<n-1$ where $n=\dim X_{\bR}$.

\medskip
(ii) In the case of exceptional type, the $D'_j$ for $j\in J_{\bR}$
are all exceptional divisors.

\medskip
(iii) In the case of ordinary type, we have
$\cJ(X_{\bR},f^{\alpha})\subset f\cO_{X_{\bR}}$ for $\alpha\ge 1$,
and hence
$$
f\cJ(X_{\bR},f^{\alpha})=\cJ(X_{\bR},f^{\alpha+1})\,\,\,\hbox{for}
\,\,\,\alpha>0,
\leqno(3.2.1)
$$
shrinking $X_{\bR}$ to an open neighborhood of the points where the
dimension of $D_{\bR}$ is $n-1$.

\medskip
(iv) The above equality (3.2.1) always holds in the complex case, and
$$
\JN(f_{\bC})=(\JN(f_{\bC})\cup(0,1])+\bN,
\leqno(3.2.2)
$$
where $\JN(f_{\bC})$ is the set of jumping numbers of $f_{\bC}$.

\medskip
The following Proposition implies the first and second assertions of
Theorem~1 in the case $n>d$, since we have always $\lct(f_{\bC})\le 1$.

\medskip\noindent
{\bf 3.3.~Proposition.} {\it
If $f$ has only an isolated zero of simple type {\rm (3.1)}, then
$$
\cJ(X_{\bR},f^{\alpha})_0=\fm_0^{[\alpha d]-n+1},\quad
{\rm RJN}(f)=\{k/d\,|\,k\ge n\},\quad
\rlct(f)=n/d,
$$
where $\fm_0$ be the maximal ideal of $\cO_{X_{\bR},0}$ and
${\rm RJN}(f)$ denote the set of real jumping numbers of $f$.}

\medskip\noindent
{\it Proof.}
In this case, we get a real resolution of singularities by the
blow-up along the origin, and $D_{\bR}=\{0\}$
since the exceptional divisor is the total transform of $D_{\bR}$.
In particular, $f$ is of exceptional type, see (3.1).
Then $J_{\bR}=\{1\}$ and $(m_1,a_1)=(d,n-1)$.
So the assertion follows from Proposition~1.

\medskip\noindent
{\bf 3.4~Example.}
Assume we have an expansion
$$
\hbox{$f=f_{d_1}+\sum_{k\ge d_2}f_k\in\bR\{\!\{x_1\dots,x_n\}\!\},$}
$$
with $f_k$ homogeneous of degree $k$, $f_{d_1}=g^e$ with $g$
irreducible, $h:=f-g^e\,(=\sum_{k\ge d_2}f_k)$ is nondivisible by
$f_{d_1}$, and $d_1<d_2$.
Let $Y\subset\bP^{n-1}$ be the projective hypersurfaces defined
by $g$.
Assume
$$
c:=d_2-d_1\ge e\ge 2, \,\,\, n>d:=d_1/e,
$$
and $Y_{\bR}$ is empty in the notation of (1.1).
Then $f$ has an isolated zero of simple type at the origin, and
$$
\rlct(f)>\lct(f_{\bC}),
$$
restricting $f$ to a sufficiently small
Zariski-open subset $X$ of the affine space $\bA^n$ containing the
origin and such that it is the only singular point of $f$.

Indeed, let $\pi:X'\to X$ be a resolution of singularities as in (1.2).
Here we may blow-up along the origin first.
Let $D'_1\subset X'$ denote the proper transform of the exceptional
divisor of this blow-up.
The pull-back of $f$ by the blow-up along the origin is locally given
by $x^{d_1}(y^e+x^cz)$, where the exceptional divisor is locally
defined by $x=0$, and the proper transforms of $g$ and
$h$ are locally given by $y$ and $z$ respectively.
So the intersection of the proper transform of $D$ and the exceptional
divisor by the blow-up along the origin is identified with $Y$, and
the total transform of $D$ is not a divisor with normal crossings
at the generic point of $Y$ since $c\ge e\ge 2$ and $h$ is
nondivisible by $f_{d_1}$.
So we have to blow-up along the proper transform of $Y$
(after making it smooth).
Let $D'_2\subset X'$ denote the proper transform of the exceptional
divisor of this blow-up.
Then we have
$$
\rlct(f)=\frac{a_1+1}{m_1}=\frac{n}{d_1}>
\frac{a_2+1}{m_2}=\frac{n+1}{d_1+e}\ge\lct(f_{\bC}).
\leqno(3.4.1)
$$
This also implies the first assertion of Theorem~1 with
$\rlct(f)<1$ if $n<d_1$.

\medskip\noindent
{\bf 3.5~Example.}
With the above notation and assumptions, assume further
$$
h=f_{d_2}, \,\,\, n=3, \,\,\, c=d=e=2,
$$
and $Y_{\bC}$ is smooth and intersects $Z_{\bC}$ at smooth
points of $Z_{\bC}$, where $Z$ is the hypersurface defined by $h$.
Then the resolution $\pi:X'\to X$ is obtained by the two blowing-ups
in Example~(3.4), and we have $J=\{1,2\}$, $J_{\bR}=\{1\}$, $m_1=4$,
$m_2=6$.
So $f$ has an isolated singularity at the origin, and the eigenvalues
$\lambda$ of the Milnor monodromy on $H^2(F_0,\bC)$ satisfy
$\lambda^4=1$ or $\lambda^6=1$, where $F_0$ denotes the Milnor fiber.

For $\lambda=i$, the $\lambda$-eigenspace of the Milnor cohomology
$H^2(F_0,\bC)_{\lambda}$ is calculated by the filtered de Rham
complex of a filtered simple regular holonomic $\cD$-module $(M,F)$
on $\bP^2_{\bC}$ whose restriction to the complement of $Y_{\bC}$ is
a complex variation of Hodge structure of type $(0,0)$ and rank 1,
and whose local monodromy around $Y_{\bC}$ is $-1$,
see [9] (or [8], 3.3 and 3.5).
Since $F_0M$ is a line bundle such that
$\otimes^2 F_0M=\cO_{\bP^2}(Y)$, we have $F_0M=\cO_{\bP^2}(1)$.
Since $\Gamma(\bP^2,\Omega^2_{\bP^2}(1))=0$, this implies
$$
F^2H^2(F_0,\bC)_{\lambda}=0.
$$
Thus $\rlct(f)\,(=3/4)$ does not appear in the spectrum [9] of
$f_{\bC}$.
Then $3/4$ is not a jumping number of $f_{\bC}$ by [3]
in the isolated singularity case.
Since the minimal jumping number of $f_{\bC}$ is $2/3$ by (3.4.1),
we get by (3.2.2)
$$
\JN(f_{\bC})\subset\Big\{\frac{k}{6}+j\,\Big|\,k=4,5,6;\,j\in\bN\Big\},
$$
(In fact, we can show the equality.)
On the other hand, we have by Proposition~(3.3)
$$
\RJN(f)=\Big\{\frac{k}{4}\,\Big|\,k\ge 3\Big\}.
$$
This implies the last assertion of Theorem~1.

\end{document}